\documentclass{amsart}
\usepackage[all]{xy}
\SelectTips{cm}{}
\usepackage{amsmath}
\usepackage{amssymb}
\usepackage{amscd}
\usepackage{amsthm}
\usepackage{amsfonts}
\usepackage{amsxtra}

\theoremstyle{plain}

\newtheorem{Thm}{Theorem}[section]

\theoremstyle{definition}

\newtheorem*{Rem-intro}{Remark}

\newcommand{\ZZ}{{\mathbb{Z}}}

\newcommand{\CC}{{\mathbb{C}}}

\newcommand{\RR}{{\mathbb{R}}}

\newcommand{\KK}{{\mathcal{K}}}

\begin{document}

\title[  A Case Study in Non-Commutative Topology \\\today ] 
{   A Case Study in Non-Commutative Topology\\ \today}


\author[Schochet]{Claude L.~Schochet}
\address{Department of Mathematics,
     Technion,
     Haifa 32000, Israel}

\email{clsmath@gmail.com}

\thanks{ It is a pleasure to thank Marc Rieffel who, besides contributing the most important theorems to this note, was very helpful 
in its preparation. }
\begin{abstract}
This is an expository note focused upon one example, the irrational rotation $C^*$-algebra.   We discuss how this algebra arises in nature - in quantum mechanics, 
group actions, and foliations, and we explain how $K$-theory is used to get information out of it. 

Our goal is to write as if we are sitting 
in Starbucks and explaining an idea to a good friend (on napkins, of course).  So we are interested in getting an idea across but not at all interested 
in the technical details that, in any event, would be lost if the coffee spilled.  So come with us for a drink at Starbucks!

\end{abstract}
\maketitle
\tableofcontents
\section{Introduction}

This is the opposite of a survey paper.  Here we are interested in one example, usually known as   {\emph{the irrational rotation 
$C^*$-algebra}} and written $A_\lambda $ where $\lambda $ is some real (usually irrational) number between $0$ and $1$. 
 There are a lot of choices for $\lambda $, so saying {\it {the}} irrational rotation $C^*$-algebra 
is already misleading. We will  exaggerate and simplify a lot in this paper.  So don't rely on this paper for 
theorems.  Think of it instead as your friend sitting across the table at Starbucks, giving you the right idea but .... leaving out 
the details.\footnote{We learned this technique from Dror Bar-Natan, who gave a great talk entitled  "{\emph{From 
Stonehedge to Witten, Skipping all the Details}}".} 

We will show that the irrational rotation $C^*$-algebra arises in three quite different contexts (there are more as well):
\begin{enumerate}
\item quantum mechanics
\item  action of a group on a compact space
\item foliations
\end{enumerate}

Then we will show how to get some information out of the irrational rotation $C^*$-algebra - in fact we will sketch (mostly) 
how to retrieve $\lambda $ using $K$-theory for $C^*$-algebras, which means we are almost preserving $A_\lambda $ up to isomorphism. 
Thus it turns out that as $\lambda $ varies among the irrationals between $0$ and $1/2$ there are uncountably many isomorphism 
classes of the algebras $A_\lambda $. 

This paper, then, is an advertisement for $K$-theory   by showing exactly one application of it.

\section{Quantum Mechanics}

In  1926-7 the quantum-mechanical revolution in physics changed 
our understanding of the world. As has been the pattern since, the physicists knew what they wanted, and the mathematicians
 struggled to keep up, to keep the physics honest (as a mathematician would put it). 

The simplest model of the hydrogen atom 
revolved about two operators $P$ and $Q$ that were to measure   position and momentum of the electron. 
Heisenberg and Max Born showed that if Q is the position operator  
  and P the momentum operator, then we have the {\it{ canonical commutation relation } }
\[
  PQ - QP = -i\hbar
\]       
where $\hbar $  is Planck's constant.

That is, not only is $PQ$ unequal to $QP$  (the order that you measure things makes a difference) but the difference was governed 
by a precise formula.  
(Note, by the way, that asking a physicist whether Planck's constant is rational or irrational  will get you a look of incredulity to this day. 
And this doesn't even take into account the folks who like to \lq\lq let Planck's constant go to zero".)\footnote{But physicists
do the latter all the time. It is Bohr's \lq\lq correspondence principle\rq\rq or passage to the semi-classical limit, \lq\lq semi\rq\rq\,
   because one keeps the Poisson bracket that is the shadow of the operator commutant, so physicists are quite comfortable 
with all this!}

Since in those days matrices were far more popular than linear operators, a lot of people tried to find matrices $P$ and $Q$ that satisfied the Heisenberg equation. 
If you try your hand with $2 \times 2$ matrices, for instance, you will see that this is not such an easy problem.  This came to an abrupt 
halt in 1927 when  Weyl and von Neumann
 observed (cf. \cite{vN},  \cite{W}) that it was impossible to find finite-dimensional matrices that do the job. The argument 
is very simple.   Suppose that $P$ and $Q$ are $n \times  n$ matrices with $PQ - QP = \lambda I$ (where $I$ is the identity matrix).  Take the trace of both sides (simply 
add up the elements on the main diagonal of the matrices)   to obtain 
\[
Trace(PQ) - Trace (QP) = Trace (\lambda I) = n\lambda .
\] 
However, $Trace(PQ) = Trace (QP) $ for any finite-dimensional matrices, and so $n\lambda = 0$, implying that $\lambda = 0$ and $PQ = QP$. 

We conclude  that one must use infinite-dimensional matrices. 
So the better thing to do is to take $P$ and $Q$ to be   self-adjoint unbounded operators on infinite-dimensional complex 
Hilbert space.   Following  Weyl \cite{W} , we set 
\[
U_s = exp(isP)  \qquad  and  \qquad V_t  = exp(itQ)  .
\]
The Stone - von Neumann theorem (cf. \cite{Rosenberg})  tells us that all such pairs of one-parameter 
unitary groups  are unique up to unitary equivalence.
On setting $s = t=  1$ and $\hbar  = \lambda $ we obtain 
\[
UV   = e^{2\pi i \lambda }VU
\]
This is called the Weyl form of the canonical commutation relation.
These operators are  bounded unitary operators on the same Hilbert space, $U, V \in \mathcal B(\mathcal H)$.

So we may take the (non-commuting) polynomial algebra generated by $U$, $V$, and their adjoints. We then close 
up this algebra with respect to the operator norm and reach our goal, the $C^*$-algebra $A_\lambda $, constructed visibly 
as a norm-closed, $*$-closed subalgebra of $B(\mathcal H)$.

This is the first construction of the $A_\lambda $.   We may restrict attension to $\lambda \in [0,1) $ and ask 
an elementary question: as $\lambda$   changes, how is $A_\lambda $ affected?  It turns out that the case 
of greatest interest is when $\lambda $ is irrational, and so we will restrict to that case as needed.

\section{Homeomorphisms of the Circle}
Let  $\phi : S^1 \to S^1$ be rotation of the circle  by $2\pi \lambda $ radians 
counterclockwise. Any rotation is a homeomorphism, and thus determines an action of the integers on the circle by 
sending $n$ to $\phi ^n$.  This defines an action of the integers on  $C(S^1)$, continuous complex-valued functions on the 
circle, and from this we will construct a $C^*$-algebra 
\[
C(S^1) \rtimes \ZZ  
\]
as follows.

Take $\mathcal H$ to be the Hilbert space $L^2(S^1)$ 
and let $T  \in \mathcal B(\mathcal H)$ be the bounded invertible operator corresponding to   rotation by $\phi $.  Any $f \in C(S^1)$ gives a pointwise multiplication operator
 $M_f  \in \mathcal B(\mathcal H)$.  Then $C(S^1) \rtimes \ZZ  $ is the norm-closed $*$-algebra generated by $T$ and by all of the $M_f$.  Note that 
finite sums of the form 
\[
\Sigma _{n = -k}^k  M_{f_n}T^n 
\]
are dense in $C(S^1) \rtimes \ZZ  $.  There is a unique normalized   trace\footnote{A {\emph{trace}} is an linear 
functional on the positive elements of the $C^*$-algebra taking values in $[0, +\infty ]$   and 
satisfying $\tau (x^*x) = \tau(xx^*)$ for all $x$ in the $C^*$-algebra.}
 $\tau $ on $C(S^1) \rtimes \ZZ   $
given on finite sums by
\[
\tau (\Sigma _n M_{f_n}T^n ) = \int_{S^1} f_0(t)dt  \,\,\in \RR
\]
where $dt$ is normalized Lebesgue measure on the circle. 
It is not at all hard to prove that (for $\lambda $ irrational, for which the action of $\ZZ $  on the circle is free)
\[
A_\lambda   \cong  C(S^1) \rtimes \ZZ 
\]
and in fact they have the same universal property. 

\section{Foliated Spaces} 

The local picture of a foliated space is $\RR ^p \times N$, where $N$ is some topological space. A subset of the form $\RR^p \times {n}$ is called a {\emph{plaque}}
and a measurable subset $T$ which meets each plaque   at most countably often (the simplest being $\{x\} \times N$ is called a {\emph{transversal}}.  The global 
picture is more complicated.  We say that  a (typically compact) space $X$ is a {\emph{foliated space}} if each point in $X$ has an open 
neighborhood homeomorphicm to the local picture and if locally the plaques fit together smoothly. A {\emph{leaf}} is a maximal union of overlapping plaques and it
by construction is a smooth $p$-dimensional manifold.  See \cite{MS} for details and lots of examples.

{\bf{NOTE  :   
I AM ATTACHING   PICTURES OF THE STARBUCKS MUGS THAT I GOT FROM THEIR WEBSITE THAT I WOULD LIKE TO USE FOR AUTHENTICITY!!!  IF THE
EDITORS OR STARBUCKS VETOES THIS THEN WE WILL HAVE TO GET SOMEBODY TO MAKE REAL PICTURES OF FOLIATIONS FOR US - I DON'T KNOW 
HOW TO DO THIS.
}}

Here is an example. Take a cylinder, which we shall visualize as the outside of a coffee mug with the bottom and handle 
removed!! 	We  can foliate this in several ways. 
\begin{enumerate}
\item  Take this tall mug (INSERT PICTURE OF TALL MUG)) and cut out the bottom and throw away the top.  That leaves a cylinder.
We can think of the cylinder being made up of circles- the circle that you see at each end and the infinite number of circles you would get by cutting through with a saw 
parallel to the circular end.  In this case each leaf is diffeomorphic to a circle.
\item We  can think of the cylinder as being made up of straight lines- the lines that start at one end of the cylinder and extend perpendicularly to the other end. ( as you would see on   what's left of this Starbucks mug  (INSERT PICTURE OF RED MUG)  after removing the handle and the bottom.) 
\end{enumerate}
More precisely, we can start with the space $[0,2\pi ] \times [0,1] $ sitting in the first quadrant, and then define $X$ to be the quotient space obtained by identifying the point $(x, 0)$
to the point $(x,1)$.   Then the first foliation corresponds to taking the plaques to be of the form $(x,t)$ where $t$ varies from $0$ to $1$ and the second foliation 
corresponds to taking the plaques to be of the form $(t,y)$ where $t$ varies from $0$ to $2\pi $.  

Next step.  We want to take the cylinder and glue the left and right ends together.  Precisely, glue the point $(0,t)$ to the point $(1,t)$.  This gives us a torus (aka the crust of a doughnut) 
and if you think carefully you will see that in both examples  (1) and (2) we wind up with a foliation by circles.  In fact the two foliations are mutually perpendicular.

(3)    Now the critical step.  Instead of gluing $(0,t)$ to the point $(1,t)$ we will glue $(0,t)$ to the point $(1, t + \lambda )$ (where $"+"$ means addition mod $2\pi $. Now something 
quite unusual happens and here we must specify whether $\lambda $ is rational or not.  If $\lambda $ is rational then each leaf of the foliation is actually a circle, though 
you will have to go around the torus several times to show this.   On the other hand, if $\lambda $ is irrational, then you do not get circles: every leaf is a line $\RR $.  Furthermore, 
each line is wrapped about the torus infinitely often, so that the line is actually dense in the torus. This construction is called the {\it{Kronecker flow on the torus}} (though we  don't know 
if this is in Kronecker's honor or if he actually invented it) and may also be described in terms of the differential equation $dy = \lambda dx$. 

Every foliated space satisfying very minimal technical assumptions has a $C^*$-algebra associated to it.  This is due to H. E. Winkelnkemper and to A Connes. 
The procedure has two steps. The first step is to associate a topological 
groupoid to the foliation. Then a general construction assigns a $C^*$-algebra to the groupoid.  In our context these are always {\emph{stable}}
$C^*$-algebras; they have the form $A \otimes \KK$ where $\KK $ is the $C^*$-algebra of the compact operators.\footnote{If you are an analyst you can think of $\KK$  as the 
smallest $C^*$-algebra inside $\mathcal B (\mathcal H)$ that contains all of the operators with finite-dimensional range. If you are an algebraist at heart then you will 
be pleased to hear that $\mathcal B (\mathcal H)$ is   a local ring, and $\KK $ is its unique maximal   ideal.}  Here are some 
examples:
\begin{enumerate}
\item If $X = F \times B $ for some smooth manifold $F$ and compact space $B$
with leaves of the form $F \times \{b\}$,
 then the foliation algebra is $C(B)\otimes \KK$.
\item More generally, if  $X$ is the total space of a compact  fibre bundle  
\[
F \to X \overset{\pi}\longrightarrow B
\]
then it is a foliation, where the leaves of the foliation are the subsets of $X$ of the form $\pi ^{-1}(b)$. The foliation 
algebra  is simply $C(B) \otimes \KK $.  
\item If $X$ is the torus foliated by circles, as above, then this is a special case of (2) and hence the foliation algebra is $C(S^1)\otimes \KK$.
\item (the punch line) If $X$ is the torus foliated by lines as constructed above at irrational angle $\lambda $ then the foliation algebra is $A_\lambda \otimes \KK $.
\end{enumerate}

There is a natural trace that arises in this construction as well. What is needed is an {\emph{invariant transverse measure.}  A {\emph{transversal}}
is a measurable set that meets each leaf of the foliation at most countably many times.  A {\emph {transverse measure}} measures transversals, naturally enough. 
If it has enough nice properties then it is an {\emph{invariant transverse measure.}} Not all foliations have them, but the ones we are looking 
at do.
In the case of the fibre bundle above, the foliation 
algebra  is simply $C(B) \otimes \KK $ and invariant transverse measures correspond to certain measures on $B$. 
Invariant transverse measures correspond to   Ruelle-Sullivan currents in foliation theory (cf. \cite{MS} Ch. IV.)

 In the  case of the Kronecker flow on the torus, the invariant transverse measure may be constructed from Lebesgue measure on a transverse 
circle to the foliation.   This passes to a trace on the foliation algebra which corresponds to the trace constructed above. 

Suppose that $p \in A$ is a projection and we have normalized the trace $\tau $ so that $\tau (1) = 1$.  Then $0 \leq \tau (p) \leq 1 $ by elementary considerations.  But what is the range of the map?  In the case $A = M_n(\CC )$ the range of $\tau $ would be $\{0, \frac 1n, \frac 2n, \dots , 1\}$.  What happens for $A_\lambda $?  Stay tuned.

Lest the reader feel cheated that we are not obtaining $A_\lambda $ on the nose, let us hasten to point out that if $A\otimes \KK  \cong B \otimes \KK $
(this is called {\emph{stably isomorphic}}) then $A$ and $B$ are strongly Morita equivalent, and conversely ($A$ and $B$ being separable) by deep results 
of Brown, Green, and Rieffel \cite{BGR}. So $A_\lambda $ and $A_\lambda \otimes \KK $ have the same representation theory, essentially, and, as we shall see, the same $K$-theory.

To summarize, we have shown that the $C^*$-algebra $A_\lambda $ arises in three disparite arenas of mathematics. (There are others as well, but this should be enough 
to convince you that it happens a lot!)   At this point, though, it is not at all clear to what extent the algebra is dependent upon $\lambda $.   Let's find out. 

\section{The World's Fastest Intro to $K$-theory}

Suppose first that $A$ is a unital $C^*$-algebra.  A {\emph{projection}} $p$ is an element of $A$ that satisfies $p^2 = p = p^*$.  
There are always projections, namely $0$ and $1$.  If $X$ is a connected space then these are the only projections in $C(X)$. On the other 
hand, $M_n(\CC )$ has lots of projections: for instance, take a diagonal matrix that has only ones and zero's on the diagonal. It turns out 
that $C(X)\otimes M_n(\CC)$  can have very interesting projections -   these correspond to vector bundles over $X$. 

Let  $P_n(A)$ denote the set of projections in $A\otimes M_n(\CC)$, and define $P_\infty (A)$ to be the union of the $P_n(A)$ (where 
we put $P_n $ inside of $P_{n+1}$ by sticking it in the upper left corner and adding zeros to the right and below.)   Unitary equivalence and saying that $p$ is equivalent to $p\oplus 0$ puts 
a natural equivalence relation 
 $ \sim $  on 
$P_\infty (A) $. 
 Then $P_\infty (A)/\sim$ \,\, has 
a natural direct sum operation, and we can turn it into an abelian group by doing the so-called Grothendieck construction (taking formal differences 
of projections). If you don't like that, take the free abelian group on the equivalence classes and then divide out by the subgroup generated 
by all elements of the form $[P+ Q] - [P] - [Q] $.  This gives an abelian group denoted $K_0(A)$.\footnote{A projection $p$
in $A\otimes M_n(\CC)$ corresponds  to a finitely generated projective  $A$-module in a natural way (think of 
endomorphisms of $A\oplus \overset{n}\dots \oplus A$) and hence we are really using the classical definition of 
$K_0(A)$ for arbitrary unital rings, except for the fact that there one would use all idempotents, not just 
projections. The two definitions are 
equivalent, though, by a version of Gram-Schmitt orthogonalization.  Similarly one can use invertibles rather than 
unitaries in the homotopy definition of $K_*$.}    

We may regard $K_0$ as a functor on unital $C^*$-algebras and maps, since if $f: A \to A'$ is unital then $f$ takes projections to projections, 
unitaries to unitaries, and preserves direct sum.  If $A$ is not unital then we may form its unitization $A^+$  (for example, $C_o(X)^+ \cong C(X^+)$ where $X$ is 
locally compact and $X^+$ is its one-point compactification), and then  define $K_0(A) $ to be the kernel of the map 
\[
K_0(A^+) \longrightarrow  K_0(A^+/A)  \cong \ZZ . 
\]
Note that if $A$ is separable then there are at most countably many equivalence classes of projections, and hence 
$K_0(A)$ is a countable abelian group.

For example,  take $A = \CC $. Then $P_n(A)$ consists of all of the projections in $M_n(\CC )$.  We learned in the second semester of Linear 
Algebra   that every projection is unitarily equivalent to a diagonal matrix of the form $diag(1,1,\dots 1, 0,0 ,0)$ and hence we may 
regard the equivalence classes of $P_n(\CC )$ to be the integers $\{0,1,2,\dots , n\}$.  
Then the equivalence classes of  $P_\infty (\CC )$  are classified by the ranks of the matrices which 
correspond to all of the natural numbers $\{0,1,2,\dots \}$ and 
taking formal inverses  we obtain $K_0(\CC ) \cong \ZZ $.   Note that the same answer 
emerges if we take $A = M_j(\CC)$ for any $j$, since \lq\lq matrices of matrices are matrices."  Similarly (but this requires a little work) 
we have the same answer if $A = \KK $ the compact operators.  Actually we need something stronger that is based upon this idea, namely 
this fact:
\[
  K_0(A)   \,\cong\,  K_0(A \otimes M_n(\CC))  \,\cong\,  K_0(A\otimes \KK )
\]
which we will use without further comment.
Next, we note that for commutative unital $C^*$-algebras $A = C(X)$ with associated maximal ideal space the compact space $X$, then
\[ 
K_0(C(X)) \cong K^0(X)
\]
where $K^0(X)$ is the Grothendieck group generated by complex vector bundles over $X$.

Today we will need only $K_0$ but there is a $K_1$ as well.\footnote{If $A$ is unital then define $U_n(A) $ to be the group of unitaries in 
$A \otimes M_n(\CC )$ and for $n> 0 $ let 
\[
K_j(A) = lim _n \pi _{j-1}(U_n(A)).
\]
This gives an infinite collection of groups, but Bott showed that the groups are periodic: $K_{j+2}(A) \cong K_j(A) $ so there are really just 
two groups up to isomorphism.  Further, if $A$ is separable then each $K_j(A)$ is countable.}

\section{The $K$-theory of the irrational rotation $C^*$-algebra: the bad news}

 Now, what happens to the irrational rotation $C^*$-algebra?   A seemingly elementary question arises first: does $A_\lambda $ 
have any non-trivial projections?    This was open for several years, and it led to decisive work by  the second author whose results, 
together with those of Pimsner-Voiculescu, we now describe.\footnote{ Remember that it is easy to show that $A_\lambda \otimes\KK $ has 
projections, and those projections determine the $K$-theory. It is a much deeper problem  to deal with $A_\lambda $ itself.} 
We are altering the historical order a bit in what follows - see Rieffel \cite{R} for the truth.

If $\lambda $ is irrational then 
Pimsner and Voiculescu showed \cite{PV} that
\[
K_0(A_\lambda ) \,\cong\, \ZZ \oplus \ZZ
\]
independent of $\lambda $.   So using $K_0$ by itself we cannot distinguish the various $A_\lambda $. 

\section{Traces to the rescue}

If $A$ is any $C^*$-algebra with a nice trace  and $p$ and $q$ are orthogonal projections in $A$ then 
\[
\tau (p \oplus q) = \tau (p) + \tau (q)
\]
 and so the trace gives us a homomorphism 
\[
K_0(A) \overset{\tau }\longrightarrow \RR 
\]
of abelian groups.  We have remarked previously that if $A$ is separable (which we assume henceforth) 
then $K_0(A) $ is a countable abelian group, and hence $\tau (K_0(A)) $, the image of 
\[
\tau : K_0(A) \longrightarrow \RR 
\]
 is a countable subgroup of $\RR $.   There are {\emph{a lot}}  of countable subgroups of $\RR$ (cf. \cite{Fuchs})!

However, there is good news.
Pimsner and Voiculescu \cite{PV} showed that the image of the trace
\[
K_0(A_\lambda ) \cong \ZZ \oplus \ZZ  \overset{\tau}\longrightarrow \RR
\]
takes values  $\ZZ + \lambda \ZZ $, the subgroup 
of $\RR $ generated by $1$ and by $\lambda $. 
Now Rieffel had previously shown that every element of $(\ZZ + \lambda \ZZ) \cap [0,1]$ is in the range of a projection in $A_\lambda$. 
Combining these results gives us this omnibus isomorphism theorem:

\begin{Thm} (Rieffel \cite{R}, Pimsner-Voiculescu \cite{PV})  
\begin{enumerate}
\item If $\lambda $ is irrational then the image of 
\[
K_0(A_\lambda )    \overset{\tau}\longrightarrow \RR
\]
is exactly $\ZZ + \lambda \ZZ $.
\item There are uncountably many isomorphism classes of algebras $A_\lambda $ as 
$\lambda $ ranges among the irrational numbers in 
the interval $[0. 1/2 ]$.
\item If $\lambda$ and $\mu $ are irrational numbers in the interval $[0, \frac 12]$ and $A_\lambda \cong A_\mu $ then $\lambda = \mu$.
\footnote{For other irrational numbers $\nu $, take the fractional part $\{\nu \}$ and then use either $\{\nu \}$   or $1 - \{\nu \}$ to get 
into the $[0, \frac 12]$ range.}
\item  If $\lambda$ and $\mu $ are irrational numbers in the interval $[0,  \frac 12]$  and $m$ and $n$ are positive integers with 
\[
A_\lambda \otimes M_m(\CC ) \cong A_\mu \otimes M_n(\CC )
\]
then $\lambda = \mu $ and $m = n$. 
\item The algebras $A_\lambda $ and $A_\mu $ are strongly Morita equivalent  (that is, $A_\lambda \otimes\KK \cong A_\mu \otimes \KK $)
if and only if $\lambda $ and $\mu $ are in the same orbit of the action of $GL(2, \ZZ )$ on irrational numbers by linear fractional 
transformations.

\end{enumerate}
\end{Thm}

So we see that the $A_\lambda $ retain all of the sensitive information about the angle.  
If we think back to the origins of $A_\lambda $ this seems really astonishing:

\begin{itemize}
\item  The exact value of Planck's constant (and whether or not it is rational) really does seem to make something of a difference!
\item   One irrational rotation of a circle is really not like another irrational rotation of a circle.
\item The angle of the Kronecker flow deeply effects the geometry of the foliation.
\end{itemize}

Case studies are supposed to suggest questions for further study. We hope that this note has done so.

\end{document}